\documentclass[11pt]{article}

\usepackage{amsmath,amsfonts,thc,dsfont,graphicx,a4wide,enumerate}

\theoremstyle{change}
\theorembodyfont{\itshape}
\newtheorem{theorem}{\bf Theorem}[section]

\newtheorem{lemma}[theorem]{\bf Lemma}
\newtheorem{corollary}[theorem]{\bf Corollary}

\theorembodyfont{\rm}

\newtheorem{remark}[theorem]{\bf Remark}

\newcommand{\conver}{\mathop{\longrightarrow}}
\newcommand{\inprob}{\ \conver_{\cP}\ }
\newcommand{\ineltwo}{\ \conver_{L^2}\ }
\newcommand{\indist}{\longrightarrow_{\mathcal{D}}}

\newcommand{\cN}{{\mathcal N}}
\newcommand{\cP}{{\mathcal P}}
\newcommand{\cL}{{\mathcal L}}

\newcommand{\cF}{{\mathcal F}}

\newcommand{\GN}{{\mathds{N}}}
\newcommand{\GR}{{\mathds{R}}}

\newcommand{\GZ}{{\mathds{Z}}}

\newcommand{\GD}{\mathds{D}}

\newcommand{\bE}{{\mathbb{E}}}
\newcommand{\bP}{{\mathbb{P}}}

\newcommand{\bV}{\mathbb{V}}

\newcommand{\bproof}{\noindent{\sc Proof.}}
\newcommand{\eproof}{$\Box$}

\begin{document}
	\title{Truncated moments of perpetuities\\ and a new
		central limit theorem for GARCH processes\\ without Kesten's regularity}

	\author{Adam Jakubowski\footnote{E-mail: adjakubo@mat.umk.pl} and Zbigniew S. Szewczak\footnote{E-mail: zssz@mat.umk.pl} \\[3mm]
		Nicolaus Copernicus University, Toru\'n, Poland
	}

	\date{}
	
	\maketitle
	
	\begin{abstract}
		We consider a class of perpetuities which admit direct
	characterization of asymptotics of the key truncated moment. The
	class contains perpetuities without polynomial decay of tail
	probabilities and thus not satisfying Kesten's theorem. We show how
	to apply this result in deriving a new weak law of large numbers for solutions to stochastic recurrence equations and a new central limit theorem for
	GARCH(1,1) processes in the critical case.
	\end{abstract}


	\noindent{\em MSClassification 2010:} 
	60 F 05; 60 F 17; 60 E 07; 60 J 05; 60 J 35.

\section{Truncated moments of perpetuities}

Let $(A_j, B_j),\ j=1,2,\ldots$ be a sequence of independent, identically distributed
random vectors with non-negative coordinates: $A_j \geq 0, B_j \geq 0$. Suppose that the series
\begin{equation}\label{aj2eq1}
U_{\infty} = \sum_{k=1}^{\infty} B_k \prod_{j=1}^{k-1} A_j
\end{equation}
is almost surely finite (by definition $\prod_{\emptyset} \equiv 1$). Then $U_{\infty}$
is called perpetuity for it admits a natural
interpretation in insurance and finance. We refer to \cite{EKM97} for an excellent primer on
perpetuities and to the articles \cite{GoMa00} and \cite{AIR09} for an in-depth discussion of
existence, uniqueness and related properties.
If $U_{\infty}$ exists, then its law is a distributional solution to the equation
\begin{equation}\label{aj2eq2}
U =_{\mathcal D}  A U + B ,
\end{equation}
where $U$ and $(A,B)$ are copies of $U_{\infty}$ and $(A_1,B_1)$,
respectively, and $U$ and $(A,B)$ are independent.
See \cite{BDM16} for an extensive treatment of this stochastic recurrence equation. 

We will assume the following non-degeneracy conditions:
\begin{eqnarray}\label{aj2eq3}
&&\bP\big( A = 1 \big) < 1,\\
&&\bP\big( B = 0\big) < 1.\label{aj2eq4}
\end{eqnarray}
We will also assume that
there exists a constant $\kappa > 0$ such that
\begin{eqnarray}\label{aj2eq5}
&&\bE A^{\kappa} = 1,\\
&&\bE B^{\kappa} < +\infty. \label{aj2eq6}
\end{eqnarray}
Relations (\ref{aj2eq3}) and (\ref{aj2eq5}) imply that
\begin{equation}\label{aj2eq7}
\bE A^{p} < 1,
\end{equation}
whenever $p > 0$, $p < \kappa$.
The function $\psi(p) = \bE A^p I (A > 0)$ is strictly convex in $(0,\kappa)$, $\psi(\kappa) = 1$ and
$\psi(p) < 1$ in $(0,\kappa)$. Hence
\begin{equation}\label{aj2eq9}
\bE A^{\kappa}\ln A \in (0,+\infty].
\end{equation}
Moreover, for $p\in(0,\kappa)$ (\ref{aj2eq6}) and (\ref{aj2eq7}) imply convergence
in $L^{p}$ of the series defining the perpetuity and therefore  $U_{\infty}$ is a distributional solution to (\ref{aj2eq2}). This solution is unique by
\cite{Ver79}, since $\bP\big(A > 1\big) > 0$, $\bP\big(A < 1\big) > 0$ and $B \geq 0$, $\bP\big( B = 0\big) <1$,
imply that there is no  $c\in \GR$ such that $B = c(1-A)$. In
particular, the solution satisfies
\begin{equation}\label{aj2eq6a}
\bE U^p < +\infty,\quad 0 < p < \kappa.
\end{equation}

It is well-known that condition (\ref{aj2eq5}) is crucial for
power-like behavior of tail probabilities of perpetuities. If $\ln A$
conditioned on $\{A > 0\}$ has a non-arithmetic distribution and 
\[
\bE A^{\kappa} \ln^{+} A < +\infty,
\]
then
\begin{equation}\label{aj2eq8}
\bP\big( U > t\big) \asymp C t^{-\kappa}, \textrm{\,as\ \ } t \to \infty.
\end{equation}
(here and in the sequel $f(t)\asymp g(t)$ means $f(t)/g(t)
\to 1$ as $t\to\infty$ and $f(\varepsilon)\asymp g(\varepsilon)$ means $f(\varepsilon)/g(\varepsilon)
\to 1$ as $\varepsilon \searrow 0$). This result essentially belongs to Kesten \cite{Kes73}.  We refer to \cite{Gol91} for a completely elaborated
proof of this fact, benefitting from a method developed by Grincevi\v{c}ius \cite{Gri75}. 
Notice that easy examples (arithmetic) show that in general Kesten's
result is not valid, i.e. (\ref{aj2eq8}) fails to hold.

But (\ref{aj2eq8}) can fail for non-arithmetic $\ln A$, as well. Kevei \cite{Kev16} explored  the case
\[ \bE A^{\kappa} = 1\text{ and } \bE A^{\kappa} \ln^{+} A = +\infty,\] and his main assumption was
\begin{equation}\label{eq:kevei}
\overline{H}_A(x):= \bE A^{\kappa} I\big\{ \ln A > x\big\} = \ell_0(x) x^{-\alpha},
\end{equation}
where $\alpha \in (0,1]$ and $\ell_0(x)$ is a slowly varying function. Under the extra assumption that $\bE |B|^{\nu} < +\infty$, for some $\nu > \kappa$, and a highly technical condition related to the strong renewal theorem, Kevei \cite[Theorem 1.1]{Kev16} proved that it is possible to obtain regularly varying tails of the form 
\[ \bP \big( U > t \big) \asymp D \frac{t^{-\kappa}}{m(\ln t)},\]
where $m(x)=\int_0^x \overline{H}_A(s)\,ds$. 

In \cite{Kev17} Kevei extended the results of \cite{Gri75} and \cite{JeO-C12} and gave the corresponding theory for the case of  arithmetic distribution of $\ln A$ conditioned on $\{A > 0\}$.
As expected, within such a framework the tails of the generated perpetuity {\em are not} regularly varying.  For refinements in other directions we refer to \cite{DaKo19}.

In view of the above discussion it is  interesting that the truncated $\kappa$-th moment of $U$ exhibits remarkable
regularity under minimal conditions.

\begin{theorem}\label{ThMain}
	Let $(A,B)$ be a random vector with nonnegative components satisfying
	conditions (\ref{aj2eq3}) -- (\ref{aj2eq6}).
	
	Let $U$ represent the unique distributional solution to equation (\ref{aj2eq2})
	(with $U$ and $(A,B)$ independent).
	
	Suppose that
	\begin{equation}\label{aj2eq12}
	\bE A^{\kappa} \ln^{+} \big( A \wedge t\big) = h_A(\ln t),
	\end{equation}
	where 
	\[ h_A(x) = x^{\rho} \ell(x),\]
	$0\leq \rho <1$ and $\ell(x)$ is a slowly varying function. Then, as $t\to\infty$,
	\begin{equation}\label{aj2eq13a}
	\bE U^{\kappa} I \{ U \leq t\} \asymp
	\bE \big((AU + B)^{\kappa} - (AU)^{\kappa}\big) g_A(t),
	\end{equation}
	where
	\begin{equation}\label{aj2eq13b} g_A(t) = 
	\begin{cases}
	\frac{\displaystyle\ln t}{\displaystyle\ell(\ln t)},&\text{ if } \rho = 0;\\
	\frac{\displaystyle\sin (\pi\rho)}{\displaystyle\pi \rho (1 - \rho)} \frac{\displaystyle\big(\ln t\big)^{1-\rho}}{\displaystyle\ell\big(\ln t\big)},&\text{ if } \rho \in (0,1). 
	\end{cases}
	\end{equation}

	In particular, if	
	\begin{equation}\label{aj2eq10}
	\bE A^{\kappa} \ln^{+} A < +\infty,
	\end{equation}
	then, as $t\to\infty$,
	\begin{equation}\label{aj2eq11}
	\bE U^{\kappa} I \{ U \leq t\} \asymp
	\frac{E ((AU + B)^{\kappa} - (AU)^{\kappa})}{E A^{\kappa}\ln  A} \ln  t. 
	\end{equation}
\end{theorem}

The proof is based on multiple application of the Karamata Tauberian Theorem. We postpone it till the end of this section.

\begin{remark} It is easy to check that (\ref{aj2eq8}) implies
	\[ \bE U^{\kappa} I \{ U \leq t\} \asymp \kappa C \ln  t.\]
	It follows that using (\ref{aj2eq11}) we can identify the well-known constant in (\ref{aj2eq8}):
	\[ C = \frac{\bE ((A + BU)^{\kappa} - (BU)^{\kappa})}{\kappa \bE B^{\kappa}\ln B}.\]
\end{remark}

\begin{remark} Now suppose that 
	\begin{equation}\label{eq:tailrv}
	\bP\big( U > t\big) = t^{-\kappa} \ell_1(t),
	\end{equation}
	for some slowly varying $\ell_1$ (e.g. as in \cite{Kev16}). Then by \cite[Proposition 1.5.9a]{BGT89} 
	\[ \ell_2(t) = \int_0^t \frac{\ell_1(s)}{s}\,ds\]
	is slowly varying and $\ell_1(t)/\ell_2(t) \to 0$. Therefore 
	\[ \bE U^{\kappa} I\{ U \leq t\} \asymp \kappa \int_0^t s^{\kappa - 1} \bP\big( U > s\big)\,ds = \int_0^t \frac{\ell_1(s)}{s}\,ds.\]
	It follows that we are able to identify the asymptotics of  $\ell_1(x)$ (up to equivalence).	
\end{remark}

\begin{remark} Let us consider Kevei's assumption (\ref{eq:kevei}):
	\[\bE A^{\kappa} I\big\{ \ln A > x\big\} = \ell_0(x) x^{-\alpha},\]
	where $\alpha \in (0,1]$.
	By the direct part of the Karamata Theorem
	\[ \bE A^{\kappa} \ln^+\big( A \wedge e^x \big) =  \int_0^x \bE A^{\kappa} I\{\ln A > v\}\,dv \asymp \ell_0(x) x^{1-\alpha}/(1-\alpha),\]
	if $\alpha \in (0,1)$ and 
	\[ \bE A^{\kappa} \ln^+\big( A \wedge e^x \big) \asymp \int_1^x \ell_0(v) \,dv/v\]
	is slowly varying, if $\alpha =1$.
\end{remark}

\noindent{\sc Proof of Theorem \ref{ThMain}}
First we shall establish the relation
	\begin{equation}\label{aj2eq14}
	\begin{split}
	\lim_{\varepsilon \searrow 0} \bE \Big[\big(A U + B\big)^{\kappa - \varepsilon} -
	&\big(A U\big)^{\kappa - \varepsilon}\Big]
	= \\ & =\bE\left[\big(A U + B\big)^{\kappa} - \big(A U\big)^{\kappa}\right] < +\infty.
	\end{split}
	\end{equation}
	If $\kappa \leq 1$, then $0 \leq \big(A U + B\big)^{\kappa - \varepsilon} - \big(A U\big)^{\kappa - \varepsilon}
	\leq
	B^{\kappa - \varepsilon} \leq 1 + B^{\kappa} $ and therefore 
	(\ref{aj2eq14}) holds. Now assume that $\kappa > 1$ and set 
	\[ k_0 =\begin{cases} \lfloor \kappa \rfloor &\text{ if } \kappa \not\in \GN \\
	\kappa - 1& \text{ if } \kappa \in \GN,
	\end{cases}
	\]
	so that $\kappa - k_0 > 0$ and $\kappa - \varepsilon - k_0 <1$ 
	for $\varepsilon > 0$. Then for small $\varepsilon$
	\begin{align*}
	0 &\leq \big(A U + B\big)^{\kappa - \varepsilon} - \big(A U\big)^{\kappa - \varepsilon}\\ 
	&= \big( AU + B\big)^{\kappa - \varepsilon - k_0} \big( AU + B\big)^{k_0}  - \big(A U\big)^{\kappa - \varepsilon - k_0} \big(A U\big)^{k_0}\\
	&\leq \big( AU\big)^{\kappa - \varepsilon - k_0} \big( AU + B\big)^{k_0} + B^{\kappa - \varepsilon - k_0} \big( AU + B\big)^{k_0} - \big(A U\big)^{\kappa - \varepsilon - k_0} \big(A U\big)^{k_0} \\
	&=  \big( AU\big)^{\kappa - \varepsilon - k_0}\sum_{j=0}^{k_0-1} {k_0 \choose j} \big(A U\big)^j B^{k_0 - j} + B^{\kappa - \varepsilon - k_0} \sum_{j=0}^{k_0} {k_0 \choose j} \big(A U\big)^j B^{k_0 - j} \\
	&\leq \big(1 +  \big(AU\big)^{\kappa - k_0}\big)\sum_{j=0}^{k_0-1} {k_0 \choose j}\big(A U\big)^j B^{k_0 - j} + \big(1 + B^{\kappa - k_0}\big) \sum_{j=0}^{k_0} {k_0 \choose j}\big(A U\big)^j B^{k_0 - j}.
	\end{align*}
	The reader may verify that the finite sum in the last line above is integrable by (\ref{aj2eq5}), (\ref{aj2eq6}) and  (\ref{aj2eq6a}). It follows that (\ref{aj2eq14}) is valid also for $\kappa > 1$.
	
	Let us now denote $G(\varepsilon) = 1 - \bE A^{\kappa - \varepsilon} = \bE \left(A^{\kappa} - A^{\kappa -
		\varepsilon}\right)$ and $H(\varepsilon) = \bE U^{\kappa - \varepsilon}$. Then we have for small $\varepsilon > 0$
	\begin{eqnarray*}
		G(\varepsilon) H(\varepsilon) = (1 - \bE A^{\kappa - \varepsilon}) \bE U^{\kappa - \varepsilon} &=&
		\bE \big(AU + B\big)^{\kappa - \varepsilon} - \bE \big(A U\big)^{\kappa - \varepsilon} \\
		&=& \bE \left[\big(A U + B\big)^{\kappa - \varepsilon} - \big(A U\big)^{\kappa - \varepsilon}\right],
	\end{eqnarray*}
	and relation (\ref{aj2eq14})
	states that 
	\begin{equation}\label{aj2eq15}
	\lim_{\varepsilon \searrow 0} G(\varepsilon) H(\varepsilon) =
	\bE\left[\big(A U + B\big)^{\kappa} - \big(A U\big)^{\kappa}\right] =: D  > 0.
	\end{equation}

	In order to examine the asymptotics of $G(\varepsilon)$ at $0$ set
	\[G_1(\varepsilon) =
	\bE \left(A^{\kappa} - A^{\kappa - \varepsilon}\right) I\{A > 1\}\]
	and $G_2(\varepsilon) = G(\varepsilon) - G_1(\varepsilon)$.
	Direct calculation shows that
	\begin{eqnarray*}
		G_2(\varepsilon) &=& \varepsilon \bE \left(A^{\kappa}
		\frac{1 - A^{-\varepsilon}}{\varepsilon}\right)I\{A \leq 1\} \\
		&\asymp& \varepsilon \bE A^{\kappa} \ln  A\ I\{ A \leq 1\} = - \varepsilon \bE A^{\kappa} \ln ^- A,
	\end{eqnarray*}
	where $ 0 \leq  \bE A^{\kappa} \ln ^- A < +\infty$.
	If (\ref{aj2eq10}) holds,  then also
	\[ G_1(\varepsilon) \asymp \varepsilon \bE A^{\kappa} \ln  A\ I\{ A > 1\} =
	\varepsilon \bE A^{\kappa}\ln ^{+} A,\]
	and we obtain that
	\begin{equation}\label{onH1}
	H(\varepsilon) \asymp \varepsilon^{-1} \frac{\bE((AU + B)^{\kappa} - (AU)^{\kappa})}{\bE A^{\kappa}\ln  A} = \varepsilon^{-1} C'.
	\end{equation}
	But $H(\varepsilon)$ is asymptotically the Laplace transform ${\cL}_R(\varepsilon)$ of a
	measure $R$ on $[0,+\infty)$ given by the formula
	\begin{equation}\label{onH2}
	R([0,x]) = \int_{[0,x]} e^{\kappa u} P_{\ln  U}(du).
	\end{equation}
	In fact,
	\[ H(\varepsilon) = \bE U^{\kappa - \varepsilon} I\{U \geq 1\} + \bE U^{\kappa - \varepsilon}
	I\{U < 1\} = H_1(\varepsilon) + H_2(\varepsilon),\]
	where
	\[ \lim_{\varepsilon\searrow 0} H_2(\varepsilon) = \bE U^{\kappa} I\{U < 1\}\]
	and
	\begin{equation}\label{onH3}
	H_1(\varepsilon) = \bE e^{-\varepsilon \ln  U} e^{\kappa \ln  U} I\{\ln  U \geq 0\} =
	\int_{[0,+\infty)} e^{-\varepsilon u}\, R(du).
	\end{equation}

	By the Karamata Tauberian Theorem (see e.g. \cite[Theorem 1.7.1]{BGT89}) relation (\ref{onH1}) is equivalent to
	\[ R([0,x]) = \bE U^{\kappa} I\{ 1 \leq U \leq e^x\} \asymp C' \cdot x,\quad x \to \infty.\]
	In other words
	\[ \bE U^{\kappa} I\{ U \leq t\} \asymp \bE U^{\kappa} I\{ 1 \leq U \leq t\} \asymp C' \ln  t,
	\]
	what gives (\ref{aj2eq11}).
	
	Passing to the general case let us assume that
	\begin{equation}\label{eq:infty}
	\bE A^{\kappa} \ln^{+} \big( A \wedge t\big) = h_A(\ln t) \to +\infty.
	\end{equation}
	
	Let us notice that
	\begin{eqnarray*}
		G_1(\varepsilon) &=& \varepsilon \bE \left(A^{\kappa}
		\frac{1 - A^{-\varepsilon}}{\varepsilon}\right)I\{A > 1\} \\
		&=& \varepsilon \bE e^{\kappa\ln  A} \left(\int_0^{\ln  A} e^{-\varepsilon v}\,dv\right) I\{\ln  A > 0\}\\
		&=&
		\varepsilon
		\int_0^{+\infty} e^{-\varepsilon v} (1-F(v))\,dv = \varepsilon \cL_{Q}(\varepsilon),
	\end{eqnarray*}
	where
	\[F(v) = \bE A^{\kappa} I\{ A\leq e^v\},\]
	and
	\begin{align*}
	Q([0,x]) &= \int_0^x (1 - F(v))\,dv =  \int_0^x \bE A^{\kappa} I\{A > e^v\}\,dv \\
	&= \int_1^{e^x} \bE A^{\kappa} I \{A > u\} \,du/u = \bE A^{\kappa} \int_1^{A \wedge e^x} \,du/u \\
	&= \bE A^{\kappa} \ln^+\big( A \wedge e^x \big) = h_A(\ln (e^x)\big) = x^{\rho} \ell(x).
	\end{align*}
	Again, by the Karamata Tauberian Theorem, 
	\[\cL_Q(\varepsilon) \asymp \Gamma(1 + \rho) \varepsilon^{-\rho}\ell(1/\varepsilon),\]
	hence $G_1(\varepsilon) \asymp \Gamma(1 + \rho)\varepsilon^{1-\rho} \ell(1/\varepsilon)$. 
	We have
	\[ \frac{G_2(\varepsilon)}{G_1(\varepsilon)} \asymp \frac{-\varepsilon\bE A^{\kappa}\ln^-A}{\Gamma(1 + \rho)\varepsilon^{1-\rho} \ell(1/\varepsilon)} = -\frac{\bE A^{\kappa}\ln^-A}{\Gamma(1 + \rho)} \varepsilon^{\rho}\frac{1}{\ell(1/\varepsilon)} \to 0,\ \text{ as $\varepsilon \searrow 0$},
	\] 
	because (\ref{eq:infty}) implies either $\rho \in (0,1)$ or $\rho = 0$ and $\ell(x) \to \infty,\ x\to\infty$. It follows that $G(\varepsilon) \asymp G_1(\varepsilon)$ and finally	
	\[ H_1(\varepsilon) \asymp \frac{D}{\Gamma(1 + \rho)} \frac{\varepsilon^{\rho-1}}{\ell(1/\varepsilon)}.\]
	Similarly as in the previous case, by the Karamata Tauberian Theorem we obtain
	\[ R([0,x]) = E U^{\kappa} I\{ 1 \leq U \leq e^x\} \asymp \frac{D}{\Gamma(1+\rho)\Gamma(2-\rho)} \frac{x^{1-\rho}}{\ell(x)} ,\quad x \to \infty,\]
	or
	\[ \bE U^{\kappa} I\{ U \leq t\} \asymp \bE U^{\kappa} I\{ 1 \leq U \leq t\} \asymp \frac{D}{\Gamma(1+\rho)\Gamma(2-\rho)} \frac{\big(\ln t)^{1-\rho}}{\ell(\ln t)}.
	\]
	This  proves (\ref{aj2eq13a}) and (\ref{aj2eq13b}), for  $\Gamma(1+\rho)\Gamma(2-\rho) =1$ if $\rho = 0$ and 
	\[ \Gamma(1+\rho)\Gamma(2-\rho) = \rho \Gamma(\rho) (1 -\rho) \Gamma(1-\rho) = \rho(1-\rho) \frac{\pi}{\sin(\pi \rho)},\]
	if $\rho \in (0,1)$. $\Box$

\section{A consequence: a weak law of large numbers for stochastic recursions}

Let $\{(A_j,B_j)\}$ be an i.i.d. sequence of random vectors distributed like $(A,B)$ that satisfies (\ref{aj2eq3})--(\ref{aj2eq6}) and let $\{U_j\}$ be a sequence given by the stochastic recursion equation
\begin{equation}\label{eq:a0} U_j = A_j U_{j-1} + B_j, \ \ j = 1,2,\ldots,
\end{equation}
where $U_0$ is independent of $\{(A_j,B_j)\}$ and distributed
according to the stationary distribution (\ref{aj2eq1}). 
Theorem \ref{ThMain} leads us to the following weak law of large numbers.

\begin{theorem}\label{ThSecond}
In assumptions and notation of Theorem \ref{ThMain} we have
\begin{equation}\label{eq:Th2}
\frac{U_1^{\kappa} + U_2^{\kappa} + \ldots + U_n^{\kappa}}{ n\,  g_A( n)} \inprob \bE \big((AU + B)^{\kappa} - (AU)^{\kappa}\big). \end{equation}
\end{theorem}

We shall obtain this theorem from a more general result that might be of independent interest.

\begin{theorem}\label{ThWeakMain}

Let $\{Y_j\}$ be a sequence of non-negative random variables with identical distribution $Y_j \sim Y$, $j=1,2,\ldots$. 
We assume that 
$ \ell(x) = \bE Y I\big( Y \leq x\big)$ is slowly varying and satisfies both
\begin{equation} \label{eq:2a}
\lim_{x\to\infty} \frac{\ell\big(x \ell(x)\big)}{\ell(x)} = 1,
\end{equation}
and
\begin{equation} \label{eq:2b}
\lim_{x\to\infty} \frac{\ell\big(x /\ln x\big)}{\ell(x)} = 1.
\end{equation}
Moreover, we assume that there are numbers  $0 \leq \eta < 1$, $h_0 > 0$ and $C_{\infty} > 0$ 
 such that for all $i,j\in\GN$  and $h \geq h_0$
 \begin{equation}\label{eq:a3}
 \bE \big( \chi_h(Y_i) \chi_h(Y_j)\big) - \bE \big( \chi_h(Y_i)\big) \bE \big( \chi_h(Y_j)\big) \leq C_{\infty}\,h^2\, \eta^{|j-i|},
 \end{equation}
 where for $h > 0$
 \[\chi_h(x) = \begin{cases} x, &\text{ if $|x| < h$};\\
 h, &\text{ if $x \geq h$};\\
 -h, &\text{ if $x \leq -h$}.
 \end{cases}\]
 
 Then 
 \begin{equation}\label{eq:a4}
 \frac{Y_1 + Y_2 + \ldots + Y_n}{n \ell(n)} \inprob 1.
 \end{equation}
\end{theorem}

Before proving both theorems let us make some comments.

\begin{remark}\label{rem:2a}
Conditions (\ref{eq:2a}) and (\ref{eq:2b}) are independent. For example, function $\exp\big( \ln x/\ln \ln x\big)$ satisfies
(\ref{eq:2a}) and does not satisfy (\ref{eq:2b}), while the function given by formula (\ref{eq:Ap4}) does not satisfy (\ref{eq:2b}) and it does  (\ref{eq:2a}).
\end{remark}

\begin{remark}\label{rem:3a}
Condition (\ref{eq:a3}) resembles the well-known $\alpha$-mixing at exponential rate. Notice, however, that it is automatically satisfied by pairwise independent sequences and, more generally, by negative quadrant dependent sequences. See \cite{Jans88} for an example of a stationary sequence of bounded random variables that is pairwise independent and is not exponentially $\alpha$-mixing.
\end{remark}

\noindent{\sc Proof of Theorem \ref{ThSecond}} \  \

Let $C_{(A,B)} = \bE \big((AU + B)^{\kappa} - (AU)^{\kappa}\big)$. By Theorem \ref{ThMain}
$\ell(x) = \bE U I\big(U \leq x\big) = C_{(A,B)} g_A(x)$ is slowly varying. The reader may directly verify that both 
(\ref{eq:2a}) and (\ref{eq:2b}) are satisfied. It remains to prove that (\ref{eq:a3}) holds.

By stationarity it is enough to estimate from above the quantity
\[ \sigma^h_j =  \bE \big( \chi_h(U^{\kappa}_j) \chi_h(U_0^{\kappa})\big) - \bE \big( \chi_h(U_j^{\kappa})\big) \bE \big( \chi_h(U_0^{\kappa})\big).\]
We will do that for $0 < \kappa \leq 1$ and $\kappa > 1$ separately, because the first case seems to be the most important one (see the next section for an application with $\kappa = 1$) and its proof is considerably simpler.  

So let us assume that $\kappa \leq 1$.
Iterating (\ref{eq:a0}) and using the independence of $U_0$ and $\{(A_j,B_j)\}$ we get 
\begin{align*}
\bE \big( \chi_h(U_j^{\kappa}) &\chi_h(U_0^{\kappa})\big) = \bE \big( \chi_h\Big( \big(U_0\prod_{i=1}^j A_i + \sum_{k=1}^j B_k \prod_{i=k+1}^j A_i\big)^{\kappa}\Big) \chi_h(U_0^{\kappa})\big)\\
&\leq \bE \big( \chi_h\Big( U_0^{\kappa}\prod_{i=1}^j A_i^{\kappa} \Big) \chi_h(U_0^{\kappa})\big) + \bE \Big( \chi_h\Big( \big(\sum_{k=1}^j B_k \prod_{i=k+1}^j A_i\big)^{\kappa}\Big) \chi_h(U_0^{\kappa})\Big) \\
&= \bE \big( \chi_h\Big( U_0^{\kappa}\prod_{i=1}^j A_i^{\kappa} \Big) \chi_h(U_0^{\kappa})\big) + \bE \Big( \chi_h\Big( \big(\sum_{k=1}^j B_k \prod_{i=1}^{k-1} A_i\big)^{\kappa}\Big)\Big) \bE \big(\chi_h(U_0^{\kappa})\big) \\
&\leq \bE \big( \chi_h\Big( U_0^{\kappa}\prod_{i=1}^j A_i^{\kappa} \Big) \chi_h(U_0^{\kappa})\big) + \bE \big( \chi_h\big( U_{\infty}^{\kappa}\big) \bE \big(\chi_h(U_0^{\kappa})\big) \\
&= \bE \big( \chi_h\Big( U_0^{\kappa}\prod_{i=1}^j A_i^{\kappa} \Big) \chi_h(U_0^{\kappa})\big) +  \bE \big(\chi_h(U_0^{\kappa})\big)^2.
\end{align*}
Therefore 
\[
\sigma^h_j  \leq \bE \big( \chi_h\Big( U_0^{\kappa}\prod_{i=1}^j A_i^{\kappa} \Big)\chi_h(U_0^{\kappa})\big).\]

Let us notice that for $x \geq 0$ and $\epsilon \in (0,1)$
\[ \chi_h(x) = h \cdot \chi_1(x/h) \leq h (x/h)^{\epsilon} = h^{1-\epsilon} x^{\epsilon}.\]
It follows that for $\epsilon \in (0,1)$ and $h \geq 1 = : h_0$ we have
\[ \bE \big( \chi_h\Big( U_0^{\kappa}\prod_{i=1}^j A_i^{\kappa} \Big)\chi_h(U_0^{\kappa})\big) \leq h^{2-\epsilon} \bE U_0^{\kappa\epsilon}  \Big(\bE A_1^{\kappa\epsilon}\Big)^j \leq C_{\infty} h^2 \eta^j,\]
where $C_{\infty} = \bE U_0^{\kappa\epsilon} < +\infty$ by (\ref{aj2eq6a}) and $\eta  = \bE A_1^{\kappa\epsilon} < 1$ by (\ref{aj2eq7}). 

Let us assume now that $\kappa > 1$. Then the function $\GR^+ \ni x \mapsto \chi_h(x^{\kappa})$ is a Lipschitz function with the Lipschitz constant $L_{\kappa} = \kappa h^{(\kappa - 1)/\kappa}$. Applying \cite[Proposition D.0.1]{BDM16} we obtain  
\[ |\sigma_j^h| \leq c \eta^j h L_{\kappa} = C_{\infty} h^{2 - 1/\kappa} \eta^j,\]
where $c > 0$ and $\eta = \bE A < 1$ by  (\ref{aj2eq7}). This completes the proof of Theorem \ref{ThSecond}.
\quad$\Box$\\

\noindent{\sc Proof of Theorem \ref{ThWeakMain}} \  \

 By Theorem \ref{th:A}  we have
 \[ \ell(x) = \bE Y I( Y \leq x) \asymp \ell_1(x) =  \bE \chi_x (Y),\]
 so it is enough to prove 2.7 with $\ell_1(x)$ in place of $\ell(x)$. 
 Moreover, it is easy to see that $\ell_1(x)$ satisfies both (\ref{eq:2a}) and (\ref{eq:2b}). Set $b_n =  n\ell_1(n)$.

By (\ref{eq:2a}), (\ref{eq:2b}) and (\ref{eq:Ap3}) we have:
\begin{align*}
 \ell_1\big(b_n/\ln n\big)&=\ell_1\big(n\ell_1(n)/ \ln n \big) =  \ell_1\Big(\frac{n}{\ln n}\ell_1\big(n/\ln n\big) \frac{\ell_1(n)}{\ell_1\big(n/\ln n\big)}\Big)\\
  &\asymp
\ell_1\Big(\frac{n\ell_1\big(n/\ln n\big)}{\ln n}\Big) \asymp \ell_1(n/\ln n) \\
 &\asymp \ell_1(n) \asymp \ell_1(n\ell_1(n)) = \ell_1(b_n).
\end{align*} 
This implies that for every $t > 0$
\[ \frac{\ell_1\big(t b_n / \ln n\big)}{\ell_1(b_n)} \asymp \frac{\ell_1\big(b_n/\ln n\big)}{\ell_1(b_n)}
\asymp 1.\]
Therefore there exists a sequence $a_n \searrow 0$ such that
\begin{equation}\label{eq:bb1}
\frac{\ell_1\big(a_n b_n / \ln n\big)}{\ell_1(b_n)} \to 1, \text{ as $n\to\infty$}.\end{equation}
For the sake of clarity, let us denote
   \[ Y^{(h)} = \chi_h(Y).\]
\begin{lemma}\label{lem:bb1}
	\[\frac{Y_1 + Y_2 + \ldots + Y_n}{b_n} \inprob 1\]
	if, and only if,
\[ \frac{Y_1^{((a_n/\ln n) b_n)} + Y_2^{((a_n/\ln n) b_n)} \ldots + Y_n^{((a_n/\ln n)b_n)}}{b_n} \inprob 1.
	\]
\end{lemma}
\bproof
By Corollary \ref{cor:A} $n \bP \big( Y > b_n\big) \to 0$, hence 	
\begin{align*}
\bP\Big( Y_1 + Y_2 + \ldots + Y_n &\neq 	Y_1^{(b_n)} + Y_2^{(b_n)} \ldots + Y_n^{(b_n)}\Big) \\
&= \bP \Big(\bigcup_{j=1}^n \{ Y_j > b_n\big\}\Big) \leq n P\big( Y > b_n\big) \to 0.
\end{align*}
Next let us consider
\[D_n = 
\frac{Y_1^{(b_n)} + \ldots + Y_n^{(b_n)} - Y_1^{((a_n/\ln n) b_n)} - \ldots - Y_n^{((a_n/\ln n) b_n)}}{b_n} \geq 0.\]
We have by (\ref{eq:bb1})
\[ \bE D_n = \frac{n}{b_n}\big(\ell(b_n) - \ell((a_n/\ln n) b_n)\big) = \frac{n \ell(b_n)}{b_n}\big(1 - \frac{\ell((a_n/\ln n) b_n)}{\ell(b_n)}\big) \to 0.\ \Box\]

Let us denote
\[ T_n = \frac{Y_1^{((a_n/\ln n) b_n)} + Y_2^{((a_n/\ln n) b_n)} +\ldots + Y_n^{((a_n/\ln n) b_n)}}{b_n}. \]
We have
\[ \bE T_n = \frac{n \ell_1\big((a_n/\ln n) b_n\big)}{b_n} = \frac{\ell_1\big((a_n/\ln n) b_n\big)}{\ell_1\big(b_n\big)} \frac{n \ell_1\big(b_n\big)}{b_n} \to 1.\]
It follows that we shall complete the proof of Theorem \ref{ThWeakMain} by showing that $T_n - \bE T_n \ineltwo 0$. Let $K > 0$ be such that 
\[ K\ln \eta < -1,\]
 and let \[m_n = \lceil K \ln n\rceil.\]
 We shall split the components in $\bE \big(T_n - \bE T_n\big)^2$ into two groups.
\begin{align*}
\sum_{1 \leq i,j \leq n\atop |i-j| > m_n} 
& \bE \Big(\frac{Y_i^{(\frac{a_n}{\ln n} b_n)} - \bE Y_i^{(\frac{a_n}{\ln n} b_n)}}{b_n} \Big)\Big(\frac{Y_j^{(\frac{a_n}{\ln n} b_n)} - \bE Y_j^{(\frac{a_n}{\ln n} b_n)}}{b_n} \Big) \\
&\leq C_{\infty} \frac{\frac{a^2_n}{\ln^2 n} b_n^2}{b_n^2} \sum_{1 \leq i,j \leq n\atop |i-j| > m_n} \eta^{|i-j|} \leq \frac{2 C_{\infty} }{1 -\eta} \big(n-m_n) \eta^{m_n}\\
& \leq \frac{2 C_{\infty}}{1 -\eta} \exp \big( \ln n + (\ln \eta)  K\,\ln n\big) \to 0.
\end{align*}
So we have to consider the remaining covariances only.
\begin{align*}
&\Big|\sum_{1 \leq i,j \leq n\atop |i-j| \leq  m_n} \bE \Big(\frac{Y_i^{(\frac{a_n}{\ln n} b_n)} - \bE Y_i^{(\frac{a_n}{\ln n} b_n)}}{b_n} \Big) \Big(\frac{Y_j^{(\frac{a_n}{\ln n} b_n)} - \bE Y_j^{(\frac{a_n}{\ln n} b_n)}}{b_n} \Big)\Big| \leq \\
&\leq  \sum_{i=1}^{n} \bE \Big(\frac{Y_i^{(\frac{a_n}{\ln n} b_n)} - \bE Y_i^{(\frac{a_n}{\ln n} b_n)}}{b_n} \Big)^2 \\
& + 2 \sum_{i=1}^{n-1} \sum_{j=i+1}^{(i+m_n)\wedge n} \bE \Big|\Big(\frac{Y_i^{(\frac{a_n}{\ln n} b_n)} - \bE Y_i^{(\frac{a_n}{\ln n} b_n)}}{b_n} \Big)\Big(\frac{Y_j^{(\frac{a_n}{\ln n} b_n)} - \bE Y_j^{(\frac{a_n}{\ln n} b_n)}}{b_n} \Big)\Big|\\
&\leq  n \bE \Big(\frac{Y_1^{(\frac{a_n}{\ln n} b_n)}}{b_n}\Big)^2 + \sum_{i=1}^{n-1} \sum_{j=i+1}^{(i+m_n)\wedge n}
\bE \Big(\frac{Y_i^{(\frac{a_n}{\ln n} b_n)}}{b_n}\Big)^2   + \bE \Big(\frac{Y_j^{(\frac{a_n}{\ln n} b_n)}}{b_n}\Big)^2 \\
&\leq \big(n + 2(n-1)m_n\big) \bE \Big(\frac{Y_1^{(\frac{a_n}{\ln n} b_n)}}{b_n}\Big)^2 \leq 3 \frac{a_n}{\ln n} m_n n E \frac{Y_1^{(\frac{a_n}{\ln n} b_n)}}{b_n} \\
&= 3 \Big(\frac{a_n}{\ln n} \lceil K \ln n \rceil\Big) \frac{\ell(\frac{a_n}{\ln n} b_n)}{\ell(b_n)}n \frac{\ell(b_n)}{b_n}
 \asymp 3 K a_n\to 0. 
\end{align*}
This finishes the proof of Theorem \ref{ThWeakMain}.\quad $\Box$
\begin{remark}\label{rem:4a}
It should be pointed out that in the above proof properties 
 (\ref{eq:2a}) and (\ref{eq:2b}) are used in order to cope with 
 the convergence $T_n - \bE T_n \ineltwo 0$ only. If we know more on $\{Y_j\}$ (e.g. pairwise independence) we obtain a complete analogue of the independent case, as the next theorem shows. This is not surprising, for many of results on the a.s. convergence or the convergence in probability rely on two-dimensional joint distributions only (see e.g. \cite{Matu92}, \cite[Remark 3.2, p. 276]{Gut13}).    
	\end{remark}
Recall that two random variables $X$ and $Y$ are negatively quadrant dependent (NQD) (see \cite{Lehm66}, also \cite{J-DP83}), if 
\[ \bP\big( X > x, Y > x\big) \leq \bP\big( X > x)\bP\big( Y > x\big),\quad x,y \in \GR^1.\]
 If $X$ and $Y$ are NQD, then by the well-known Hoeffding's identity
 \begin{align*}
  \bE \chi_h(X) &\chi_h(Y) - \bE \chi_h(X) \bE \chi_h(Y) \\
  &= \int_{-h}^{h} \int_{-h}^h \Big(\bP\big( X > x, Y > y\big) - \bP\big( X > x\big)\bP\big(Y > y\big)\Big)\,dx dy \leq 0,
  \end{align*}
 and therefore, if $Y_1,\ldots, Y_n$ are NQD
 \[ \bV\text{ar}\big(\chi_h(Y_1) + \ldots + \chi_h(Y_n)\big) \leq \bV\text{ar}\big( \chi_h(Y_1)\big) +  \ldots + \bV\text{ar}\big( \chi(Y_n)\big).\]

\begin{theorem}
Let $\{Y_j\}$ be a sequence of non-negative and NQD random variables with identical distribution $Y_j \sim Y$, $j=1,2,\ldots$. 
Suppose that 
$ \ell(x) = \bE Y I\big( Y \leq x\big)$ is slowly varying
and $\{b_n\}$ satisfies
\[ \frac{n\,\ell(b_n)}{b_n} \to 1, \text{ as $n\to\infty$}.\]
Then 
\begin{equation}\label{eq:axxx}
\frac{Y_1 + Y_2 + \ldots + Y_n}{b_n} \inprob 1.
\end{equation}

\end{theorem}
\bproof
Let $a_n \searrow 0$ be such that
\[ \frac{\ell(a_n b_n)}{\ell(b_n)} \to 1. \]
Then by arguments identical as in the proof of Lemma \ref{lem:bb1} (with $\{a_n\}$ alone replacing $\{a_n/\ln n\}$) convergence (\ref{eq:axxx}) holds if, and only if, 
\[
\frac{Y_1^{(a_n b_n)} + Y_2^{(a_n b_n)} \ldots + Y_n^{(a_n b_n)}}{b_n} \inprob 1.
\]
Notice that also
\[ \bE \Big(\frac{Y_1^{(a_n b_n)} + Y_2^{(a_n b_n)} \ldots + Y_n^{(a_n b_n)}}{b_n}\Big)= \frac{n\,\ell(a_nb_n)}{b_n} = \frac{n\,\ell(b_n)}{b_n}  \frac{\ell(a_nb_n)}{\ell(b_n)}  \to 1.\]
Therefore the following natural estimate completes the proof.
\[ \bV\text{ar}\Big( \frac{Y_1^{(a_n b_n)} + Y_2^{(a_n b_n)} \ldots + Y_n^{(a_n b_n)}}{b_n}\Big) \leq \frac{n \bV\text{ar}\big( Y_1^{(a_n b_n)}\big)}{b_n^2} \leq a_n \frac{ n \ell\big(a_n b_n\big)}{b_n} \to 0. \] 
\eproof

\section{Another consequence: a central limit theorem for GARCH(1,1) processes}
\label{AGP}
A sequence $\{X_j\}$ of random variables is said to be a GARCH(1,1) process if 
\begin{eqnarray}\label{eqaj1}
X_j&=& \sigma_j Z_j\;,\\
\sigma_j^2&=&\beta+\lambda X_{j-1}^2 +\delta \sigma_{j-1}^2 \;,\label{eqaj2}
\end{eqnarray}
where the constants $\beta, \lambda, \delta$ are nonnegative,
$\{Z_j\}$ is an i.i.d. multiplicative noise, $\sigma_j \geq 0$ and $X_0$ and
$\sigma_0^2$ are given and independent of $\{Z_j\}_{j\geq 1}$. If
$\delta = 0$ in (\ref{eqaj2}) then the corresponding process is
called  ARCH(1) process.

The terminology (ARCH stands for ``Autoregressive Conditionally
Hetero\-skedastic" while GARCH is the ``Generalized ARCH") was
introduced by Engle \cite{engle} and Bollerslev \cite{Bo86} in the
context of modeling volatility phenomena in econometric time
series. Engle considered only normally distributed noise variables,
but this is too restrictive and it is reasonable to assume only that
\begin{equation}\label{eqaj3}
\bE Z_j = 0, \quad \bE Z_j^2 = 1.
\end{equation}

There exists a huge literature on both theoretical and practical
aspects of GARCH processes. As an excellent  mathematical introduction to ARCH(1) processes may serve \cite{EKM97}. Mathematics of GARCH(1,1) processes is studied in detail in \cite{BoPi92}, \cite{MiSt00}, \cite{BDM02}, see also \cite[Chapters 2 and 3]{BDM16}. For financial aspects of modeling with GARCH(1,1) processes we refer to the extensive sources 
\cite{handbk} and \cite{FrZa10}.

Here we shall focus on seldom investigated properties
of GARCH processes related to the threshold condition $\lambda +
\delta =1$.

It is well known that if $\lambda + \delta < 1$, then there exists
a strictly stationary sequence $\{(X_j,\sigma_j^2)\}$ built on the
i.i.d. noise $\{Z_j\}_{j\in\GZ}$,
satisfying (\ref{eqaj1}) and (\ref{eqaj2}) and such that
\begin{equation}\label{eqaj4}
\bE \sigma_j^2 = \bE X_j^2 = \frac{\beta}{1 - \lambda - \delta}.
\end{equation}
If $\lambda + \delta > 1$ and the stationary solution exists, then it has heavy tails (see \cite{BJMW11}, \cite{BDM16}  for the corresponding limit theory with stable limits).

When $\lambda + \delta=1$ and $\beta>0$, a simple choice $Z_j = \pm 1$ with probability $1/2$ provides an example with no stationary solution. It is not difficult to show (see e.g. \cite{BDM16}) that a necessary and sufficient condition for the existence of a (unique in law) stationary distribution for 
(\ref{eqaj1})
and (\ref{eqaj2})
is that
\begin{equation}\label{eq:garch}
 \beta > 0, \text{ and } \bE\ln\big(\lambda Z_1^2 + \delta\big) < 0.
 \end{equation}
In any case, if the stationary solution exists, it is of infinite variance. This makes the modeling with GARCH a delicate problem, for estimates performed on
real data often give the value of $\lambda + \delta$ very close to $1$ (e.g. 0.995 - see \cite{CS03}, also \cite{EnBo86}). 
It follows that the critical case $\lambda + \delta = 1$ is
interesting from
the point of view of both mathematics and econometrics. 

Here we are going to prove a central limit theorem for GARCH(1,1) processes in the case
when $\lambda + \delta = 1$ and under {\em minimal assumptions} on the
marginal distribution of
the noise sequence $\{Z_j\}$.

To give a flavor of necessary reasoning we begin with
discussion of two central limit theorems for
ARCH(1) processes ($\delta = 0$). For the time being we shall assume that
\[
	\text{ the noise i.i.d. variables $\{Z_j\}$ are standard normal.}\]
Then
$\lambda <1$ implies
\begin{equation}\label{eq:clt1} \frac{X_1 + X_2 + \ldots + X_n}{\sqrt{n}} \indist
\cN\big(0,\frac{\beta}{1 - \lambda}\big),
\end{equation}
while $\lambda = 1$ implies
\begin{equation}\label{eq:clt2}  \frac{X_1 + X_2 + \ldots + X_n}{\sqrt{n\ln  n}} \indist
\mathcal{N}\big(0,C_{\beta,1}\big),
\end{equation}
 where
\[C_{\beta,1} =
\frac{\beta}{E\Big[\big(Z_1^2) \ln  (Z_1^2)\Big]}\approx 1.3705
\cdot \beta.\]
(\ref{eq:clt1}) can be proved in various ways.
One possible direction is based on mixing properties of
GARCH processes. Mikosch and St\u{a}ric\u{a} \cite{MiSt00} 
proved that GARCH(1,1) processes with Gaussian noise are strongly (or $\alpha$-)mixing with exponential rate.
This means that $\alpha(n) \leq K \eta^n$ for some constants $K > 0$ and $\eta \in [0,1)$, where
for a stochastic process $\{Y_k\}_{k\in\GN}$ the well-known coefficient $\alpha(n) = \alpha(n,\{Y_k\})$ is defined
as
\[ \alpha(n) = \sup \left\{|P(A\cap B) - P(A)P(B)| \,: \, A\in \cF_{1}^{m}, B \in \cF_{m+n}^{\infty}, 
m \in \GN \right\},\]
with $\cF_{1}^{m} = \sigma\{Y_k\, :\, k \leq m\}$ and $\cF_{m+n}^{\infty} = \{Y_k\, :\, k \geq m+n\}$ (see e.g. \cite{Brad07} or \cite{Douk94} for properties and examples).
Since we have exponential $\alpha$-mixing and there exist moments higher than $2$ (due to 
$\lambda<1$), (\ref{eq:clt1}) is a direct consequence of 
Ibragimov's CLT for strongly mixing sequences (see e.g. \cite[Theorem 18.5.3, p. 346]{IbLi71}). 

On the other hand
$\{X_{n,k} = \frac{X_k}{\sqrt{n}}\, :\,  k=1,2,\ldots,n,\, n\in\GN\}$
is a square integrable martingale difference array, so one might also use a suitable version of the Martingale CLT, as it is done later in this section.

To avoid technicalities we prefer another proof, based on the fact that the regular
conditional distribution of
$X_{n}$ with respect to the ``past" is ${\mathcal N}(0,\beta + \lambda
X^2_{n-1})$.
Let us recall a device related to the Principle of Conditioning (see \cite{AJ80}, \cite{AJ86} and \cite[Appendix]{EMJV20} for an extended version).

\begin{lemma}
	Let $\{X_{n,k}\,;\, k=1,2,\ldots,k_n,\, n\in\GN\}$ be an array of random variables which are row-wise adapted to
	a sequence of filtrations $\{\{\mathcal{F}_{n,k}\}\}_{n\in\GN}$.
	Define
	\[\phi_{n,k}(\theta) = E (e^{i\theta X_{n,k}}|\mathcal{F}_{n,k-1}),\quad
	\phi_n(\theta) = \phi_{n,1}(\theta) \cdot \phi_{n,2}(\theta)\cdot \ldots
	\cdot \phi_{n,k_n}(\theta).\] 
	If $ \phi_n(\theta) \inprob C(\theta)\neq 0$,
	then also
	\[ E e^{i\theta (X_{n,1} + X_{n,2} + \ldots + X_{n,k_n})} \to C(\theta).\]
\end{lemma}
Given the above lemma, the proof of (\ref{eq:clt1}) is in one line: setting $X_{n,k} = X_k/\sqrt{n}$ and
applying the 
individual ergodic theorem one obtains:
\[ -\ln  \phi_n(\theta) = \frac{1}{n} \sum_{k=1}^{n} \frac{1}{2}\theta^2
(\beta + \lambda X_{k-1}^2) \to \frac{1}{2}\theta^2(\beta +
\lambda E X_0^2)  =  \frac{1}{2}\theta^2\frac{\beta}{1 - \lambda} \text{ a.s. }.\]

When we try to prove (\ref{eq:clt2}) the same way, we obtain ($\lambda = 1$): 
\[ -\ln  \phi_n(\theta) = \frac{1}{n\ln  n} \sum_{k=1}^{n} \frac{1}{2}\theta^2
(\beta + X_{k-1}^2 ) \asymp \frac{1}{2}\theta^2\frac{1}{n\ln  n}\sum_{k=1}^{n} 
X_{k-1}^2\] 
and the convergence in probability of $\phi_n(\theta)$ is not obvious, unless we have at disposal a weak law of large numbers for $\{X_j^2\}$! A suitable law of large numbers and the corresponding central limit theorem (\ref{eq:clt2}) were proved in \cite[Example 1]{Szew12}.

It should be pointed out that the results of \cite{Szew12} rely heavily on the assumption of exponential $\alpha$-mixing as well as on (\ref{aj2eq8}) held for $\{X_j^2\}$ with $\kappa = 1$ (Kesten's regularity). We know from Section 1 that the power tail decay does not hold in many cases. Similarly, there seems to be no general result on exponential $\alpha$-mixing valid for all GARCH(1,1) processes (see \cite{BDM02} or \cite[Proposition 2.2.4, p. 23]{BDM16}).

The main advantage of our approach is that we can use the weak law of large numbers given in Theorem \ref{ThSecond}, where we need only natural non-degeneracy assumptions (\ref{aj2eq3})--(\ref{aj2eq6})  and a weak one-sided covariance bound given by (\ref{eq:a3}), ideally suited for stationary solutions to stochastic recurrence equations, hence also for GARCH processes. 

\begin{theorem}\label{th:CLT}
	Suppose that $\lambda > 0$, $\lambda + \delta = 1$, $\bP\big( Z_1^2 \neq 1\big) > 0$ and that
	\begin{equation}\label{cor:garcz1}
	\bE \big(1 + \lambda (Z_1^2 - 1)\big) \ln^{+} \big( \big(1 + \lambda (Z_1^2 - 1)\big) \wedge t\big) = h_A(\ln t),
	\end{equation}
	where 
	\[ h_A(x) = x^{\rho} \ell(x),\]
	$0\leq \rho <1$ and $\ell(x)$ is a slowly varying function. Then we have
\begin{align}\label{eq:CLTa}
\frac{X_1 + X_2 + \ldots + X_n}{ \sqrt{n\,  g_A(n)}} &\indist \cN\big(0,\beta\big),\\
\frac{X^2_1 + X^2_2 + \ldots + X^2_n}{ n\,  g_A(n)} &\inprob \beta,\label{eq:CLTb} \\
\frac{\sigma^2_1 + \sigma^2_2 + \ldots + \sigma^2_n}{ n\,  g_A(n)} &\inprob \beta, \label{eq:CLTc}
\end{align}
where $g_A(x)$ is given by (\ref{aj2eq13b}).

In particular, if 
\begin{equation}\label{cor:garcz2} \bE \big(1 + \lambda (Z_1^2 - 1)\big) \ln^{+}  \big(1 + \lambda (Z_1^2 - 1)\big) < +\infty,
\end{equation}
then 
\begin{align}\label{eq:CLTaa}
\frac{X_1 + X_2 + \ldots + X_n}{ \sqrt{n\,  \ln n}} &\indist \cN\big(0,\beta\, C_{\lambda, Z}\big),\\
\frac{X^2_1 + X^2_2 + \ldots + X^2_n}{ n\,  \ln n} &\inprob \beta\,C_{\lambda, Z},\label{eq:CLTbb} \\
\frac{\sigma^2_1 + \sigma^2_2 + \ldots + \sigma^2_n}{ n\, \ln n} &\inprob \beta\,C_{\lambda, Z}, \label{eq:CLTcc}
\end{align}
where
\[ C_{\lambda,Z} = \frac{1}{\bE \big(1 + \lambda (Z_1^2 - 1)\big) \ln  \big(1 + \lambda (Z_1^2 - 1)\big)}.\] 
\end{theorem}

In fact both (\ref{eq:CLTa}) and (\ref{eq:CLTaa}) can be strengthened to the functional convergence  on the Skorokhod space $\GD\big([0,+\infty)\big)$ equipped with Skorokhod's topology $J_1$. We refer to \cite{JaSh03} for necessary definitions and results.

Let us define
\[ S_n(t) = \sum_{j=1}^{\lfloor n\,t\rfloor} X_j,\quad t > 0,\]
and let $\{W(t)\,;\,t > 0\}$ denotes the standard Wiener process. 

\begin{theorem}\label{th:functional}
In assumptions and notation of Theorem \ref{th:CLT}, relation   (\ref{cor:garcz1}) implies that 
\begin{equation}\label{eq:FCLT1}
\frac{S_n(t)}{ \sqrt{n\,  g_A(n)}} \indist \sqrt{\beta}\, W(t), 
\end{equation}
and (\ref{cor:garcz2}) implies that
\begin{equation}\label{eq:FCLT2}
\frac{S_n(t)}{ \sqrt{n\,\ln  n}} \indist \sqrt{\beta C_{\lambda,Z}}\, W(t).
\end{equation}
where in both cases the convergence in law holds on the space $\Big(\GD\big([0,+\infty)\big), J_1\Big)$.
\end{theorem}
\begin{remark} Clearly, from the point of view of possible applications relations (\ref{eq:CLTaa})--(\ref{eq:CLTcc}) and (\ref{eq:FCLT2}) are the most important, for they refer to common noise variables $\{Z_j\}$ (like $\cN(0,1)$, normalized $t$-Student's distributions, etc.). On the other hand, more general relations (\ref{eq:CLTa})--(\ref{eq:CLTc}) and (\ref{eq:FCLT1}) illustrate the remarkable flexibility of the model.
	\end{remark}

\begin{remark}
With the law of $Z$ fixed, function $y(\lambda) = C_{\lambda,Z}$ is strictly decreasing from $+\infty$ (for $\lambda = 0+$) to
$1/\bE Z^2 \ln Z^2$ (for $\lambda = 1$). As an example may serve 
\[ C_{\lambda,Z} = \frac{2}{\big(1+\lambda\big) \ln \big( 1 + \lambda\big) + \big(1 - \lambda) \ln\big(1-\lambda\big)},\] 
obtained for 
\[ Z = \begin{cases} 
\sqrt{2}, &\text{ with probability $1/4$},\\
0, &\text{ with probability $1/2$},\\
		-\sqrt{2}, &\text{ with probability $1/4$}.
\end{cases}\]
Notice that in this simple example Kesten's regularity does not hold.
\end{remark}
{\sc Proofs of Theorems \ref{th:CLT} and \ref{th:functional}}

Let us consider the stochastic recurrence equation implied by (\ref{eqaj1}) and (\ref{eqaj2}) and specified for the case $\lambda > 0$, $\lambda + \delta = 1$.
\[ \sigma_j^2 = \big(\lambda Z_j^2 + \delta\big) \sigma_{j-1}^2 + \beta = \big(1 + \lambda (Z_j^2 -1)\big) \sigma_{j-1}^2 + \beta = A_j \sigma_{j-1}^2 + B_j.\]
Condition (\ref{aj2eq3}) is satisfied if $\bP \big( X_j^2 \neq 1\big) > 0$ and (\ref{aj2eq4})  holds if $\beta > 0$. Notice that these non-degeneracy assumptions exclude the trivial case when $X_j^2 \equiv 1$  and when there is no stationary solution to (\ref{eqaj1})--(\ref{eqaj2}). Condition (\ref{aj2eq5}) holds for $\kappa = 1$ and therefore
\[ \bE \big((AU + B)^{\kappa} - (AU)^{\kappa}\big) = \beta.\]
Finally (\ref{aj2eq6}) is trivial.

It follows that that we may apply Theorem \ref{ThSecond} to the sequence $\{\sigma_j^2\}$:
\begin{equation}\label{clt:e1}
\frac{\sigma^2_1 + \sigma_2^2 + \ldots + \sigma_{n}^2}{ n\,  g_A( \ln n)} \inprob \beta. \end{equation}
Let 
\[ \varsigma_{n,j} = \sqrt{\chi_1\Big(\frac{\sigma^2_{j}}{n\,g_A(\ln n)}\Big)}, \ \ n,j\in \GN.\]
We have by (\ref{clt:e1}) and Corollary \ref{cor:A}
\[ \sum_{j=1}^{n} \varsigma_{n,j}^2 \inprob \beta.\]
In fact, by the row-wise stationarity of $\{\varsigma_{n,j}^2\}$ and the $1$-regular variation of $b_n = n\,g_A(\ln n)$ we have more:
\begin{equation}\label{clt:e2}
Q_n(t) =  \sum_{j=1}^{\lfloor n\,t\rfloor} \varsigma_{n,j}^2 \inprob t \beta,\ \ t \geq 0.
\end{equation}
Set 
\[Y_{n,j} = \varsigma_{n,j} Z_j,\ \ j,n\in \GN,\quad \Sigma_n(t) = \sum_{j=1}^{\lfloor n\,t\rfloor} Y_{n,j}, \ \ t \geq 0, n\in\GN.\]
By Corollary \ref{cor:A} we have for each $T > 0$
\begin{equation}\label{almost}
 \bP\Big( \exists_{t\in [0,T]}\frac{S_n(t)}{\sqrt{n\,g_A(\ln n)}} \neq \Sigma_n(t)\Big) \leq (nT)\bP\big( \sigma_1^2 > b_n\big) \to 0,\end{equation} 
and so it is enough to prove a functional limit theorem for processes  $\Sigma_n(t)$. Notice that $\{Y_{n,j}\}$ is a martingale difference array for which  (\ref{clt:e2}) gives
the convergence of conditional variances:
\[ \sum_{j=1}^{\lfloor n t\rfloor} \bE \big( Y^2_{n,j}\big| \cF_{j-1}\big) = Q_n(t) \inprob t\beta,\ \ t\geq 0,  \]
where $\cF_{j} = \sigma\big( \sigma^2_0, Z_1, Z_2, \ldots, Z_j\big)$, $j\in\GN$.
By \cite[Theorem 3.33, p. 478]{JaSh03}  we have to check the  Lindeberg condition in the conditional form. 

Let $a_n\searrow 0$ be given by (\ref{eq:Apy}). Let us notice that
\begin{align*}
 I\big( \varsigma_{n,j}^2 Z_j^2 > \varepsilon\big) &=   I\big( \varsigma_{n,j}^2 Z_j^2 > \varepsilon, \varsigma_{n,j}^2 \leq a_n\big) +  I\big( \varsigma_{n,j}^2 Z_j^2 > \varepsilon, \varsigma_{n,j}^2 > a_n\big) \\
 &\leq  I\big(Z_j^2 > \varepsilon/a_n\big) +  I\big(\varsigma_{n,j}^2 > a_n\big).
 \end{align*}
 Therefore
\begin{align*}
\bE \big( Y^2_{n,j} I\big( Y^2_{n,j} > \varepsilon\big) \big| \cF_{j-1}\big)  &= \bE \big( \varsigma^2_{n,j}Z_j^2 I\big( \varsigma^2_{n,j} Z_j^2 > \varepsilon\big) \big| \cF_{j-1}\big) \\
& \leq \varsigma^2_{n,j} \bE Z_j^2 I\big(Z_j^2 > \varepsilon/a_n\big) + \varsigma^2_{n,j} I\big(\varsigma_{n,j}^2 > a_n\big),
\end{align*}
and
\[ \sum_{j=1}^{\lfloor n t\rfloor} \bE \big( Y^2_{n,j} I\big( Y^2_{n,j} > \varepsilon\big)\big| \cF_{j-1}\big) \leq Q_n(t)\bE Z_j^2 I\big(Z_j^2 > \varepsilon/a_n\big) + \sum_{j=1}^{\lfloor n t\rfloor}\varsigma^2_{n,j} I\big(\varsigma_{n,j}^2 > a_n\big).\]
The first term on the right hand side trivially converges in probability to $0$, while for the second term we obtain by 
(\ref{eq:Apy}) 
\[ \bP \Big( \sum_{j=1}^{\lfloor n t\rfloor}\varsigma^2_{n,j} I\big(\varsigma_{n,j}^2 > a_n\big) > 0\Big) \leq nt\, \bP \big(\sigma^2_1 > a_n b_n\big) \to 0,\ \ t > 0.\]
It follows that on the space $\Big(\GD\big([0,+\infty)\big), J_1\Big)$ \[  \Sigma_n(t) \indist \sqrt{\beta} W(t),\]
hence by (\ref{almost}) we have also (\ref{eq:FCLT1}). Applying again 
 \cite[Theorem 3.33, p. 478]{JaSh03} we obtain (\ref{eq:CLTb}).
$\Box$.

\begin{appendix}
	\section{On slowly varying functions}
In this section we gather some properties of slowly varying functions which are crucial for our reasoning.

A measurable positive function $\ell : [x_0,+\infty) \to \GR^+$, $x_0 > 0$,  is slowly varying, if for every $t > 0$
\[ \lim_{x\to\infty} \frac{\ell(t x)}{\ell(x)} = 1.\]
By \cite[Theorem 1.5.13]{BGT89} there exists the Bruin conjugate $\ell^{\#}(x)$ of $\ell(x)$ that is determined uniquely up to the asymptotic equivalence by the relations
\[ \lim_{x\to\infty} \ell(x)\ell^{\#}(x\ell(x)) =1,\ \ \lim_{x\to\infty} \ell^{\#}(x)\ell(x\ell^{\#}(x)) =1.\]
If we set $b_n = n \ell_0^{\#}(n)$, where $\ell_0(x) = 1/\ell(x)$, then by the second relation above $\{b_n\}$ satisfies 
\begin{equation}\label{eq:Ap1}
\frac{n \ell(b_n)}{b_n} \conver 1, \text{ as $n\to\infty$}.
\end{equation}
In particular, there exists a sequence $a_n \searrow 0$ such that 
\begin{equation}\label{eq:Ap2}
 \frac{\ell\big(a_n b_n\big)}{\ell\big(b_n\big)} \conver 1, \text{ as $n\to\infty$}.
\end{equation}
Indeed, $b_n \to \infty$, hence we have $\ell\big(t b_n\big)/\ell\big(b_n\big) \to 1$ for every $t > 0$. 
Therefore (\ref{eq:Ap2}) holds if $a_n \searrow 0$ slowly enough.

The construction of $b_n$ is considerably easier, if $\ell(x)$ satisfies (\ref{eq:2a}). In such a case 
\[ \frac{n \ell\big( n\ell(n)\big)}{n\ell(n)} = \frac{ \ell\big( n\ell(n)\big)}{\ell(n)} \conver 1,\text{ as $n\to\infty$},\] 
and it is enough to set 
\begin{equation}\label{eq:Ap3}
b_n = n \ell(n).
\end{equation}
It should be pointed out that not all slowly varying functions satisfy (\ref{eq:2a}). A suitable example can be taken from  \cite[p. 302]{BoSe71}:
\begin{equation}
\label{eq:Ap4}
\ell(x) = \exp\big( (\ln x)^{\beta}\big), \ \ \frac{1}{2} < \beta <1.
\end{equation}

The next fact can be deduced from \cite[p. 283, Theorem 2]{Fell71}. Since it is of crucial importance for our reasoning and since the proofs in \cite{Fell71} are a bit informal we provide here a direct proof based on core properties of slowly varying functions.

\begin{theorem}\label{th:A}
	Let $Y$ be a non-negative random variable such that
	    \[ \ell(x)= \bE Y I(Y \leq x)\]
	    is a slowly varying function.
Then
\begin{equation}\label{eq:Ap5}
\lim_{x\to\infty} \frac{x P( Y > x)}{E Y I(Y \leq x)} = 0.
\end{equation}
In particular 
\begin{equation}\label{eq:Ap6}
\ell_1(x) = E Y \wedge x = E Y I(Y \leq x) + x P(Y > x)  \asymp \ell(x).
\end{equation}
\end{theorem}
\noindent{\sc Proof.}
By the Fubini theorem we have for $x > 0$
\begin{align*} \int_x^{\infty} \big(E Y I( Y\leq y) - &E Y I( Y\leq x)\big) \frac{d\,y}{y^2}\\ 
&= E Y \int_{x\vee Y}^{\infty} \frac{d\,y}{y^2} - E Y I( Y\leq x) \frac{1}{x}\\
&= E \frac{Y}{x \vee Y} - E Y I( Y\leq x) \frac{1}{x}\\
&=  E Y I( Y\leq x) \frac{1}{x} + P ( Y > x) - E Y I( Y\leq x) \frac{1}{x} \\
&= P( Y > x).
\end{align*}
Therefore
\begin{align*} \frac{x P( Y > x)}{E Y I( Y\leq x)} &= \frac{x}{\ell(x)}  \int_x^{\infty} \big(\ell(x) - \ell(x)\big)
\frac{d\,y}{y^2}\\ 
&= \frac{x}{\ell(x)} \frac{\ell(x)}{x} \int_1^{\infty} \Big(\frac{\ell(tx)}{\ell(x)} - 1\Big) \frac{d\,t}{t^2} = 
\int_1^{\infty} \Big(\frac{\ell(tx)}{\ell(x)} - 1\Big) \frac{d\,t}{t^2}.
\end{align*}
Take $\delta \in (0,1)$. By the Potter Theorem \cite[Theorem 1.5.6]{BGT89} there exist  constants $C_1, C_2 \geq 1$ such that for $u \geq v \geq C_2$ we have
\[ \frac{\ell(u)}{\ell(v)} \leq C_1 \frac{u^{\delta}}{v^{\delta}}.\]
Take $\varepsilon > 0$ and let $T \geq C_2$ be such that 
\[ C_1\int_{T}^{\infty} \frac{d\,t}{t^{2-\delta}} < \varepsilon/2. \]
Further, by the Uniform Convergence Theorem \cite[Theorem 1.2.1]{BGT89}, let 
$x_0 \geq C_2$  be such that 
\[ \sup_{t\in [1,T]} \Big(\frac{\ell(tx)}{\ell(x)} - 1\Big) < \varepsilon/2, \ \ x \geq x_0.\]
Then for $x \geq x_0$ we have 
\begin{align*}
\int_1^{\infty} \Big(\frac{\ell(tx)}{\ell(x)} - 1\Big) \frac{d\,t}{t^2} &= \int_1^{T} \Big(\frac{\ell(tx)}{\ell(x)} - 1\Big) \frac{d\,t}{t^2} + \int_T^{\infty} \Big(\frac{\ell(tx)}{\ell(x)} - 1\Big) \frac{d\,t}{t^2} \\
&\leq 
\int_1^T \varepsilon/2\frac{d\,t}{t^2} +  \int_T^{\infty}
C_1 t^{\delta} \frac{d\,t}{t^2} = \varepsilon.
\end{align*}

\begin{corollary}\label{cor:A}
In assumptions of Theorem \ref{th:A}, if $b_n$ is given by (\ref{eq:Ap1}), then
\begin{equation}\label{eq:Apx}
n P\big( Y > b_n\big) \conver 0, \ \ \text{ as $n \to\infty$},
\end{equation}	
and there exists a sequence $a_n \searrow 0$ such that still
\begin{equation}\label{eq:Apy}
n P\big( Y > a_n b_n\big) \conver 0, \ \ \text{ as $n \to\infty$}.
\end{equation}	
	\end{corollary}

\bproof
	
We have for each $t > 0$
	\[ n P\big(Y > t\,b_n\big) = \frac{ n \ell\big(b_n\big)}{b_n}\cdot \frac{\ell(t b_n)}{t\, \ell(b_n)}\cdot \frac{t\,b_n P (Y > t\,b_n)}{\ell(t\,b_n)} \conver 0, \ \ \text{ as $n \to\infty$}.\]
	Therefore it is enough to take $a_n \searrow 0$ slowly enough.
	\eproof
	
\end{appendix}

\end{document}